\documentclass[10pt,reqno]{amsart}
\usepackage{amsmath, amsthm, amscd, amsfonts, amssymb, graphicx, color}
\usepackage[bookmarksnumbered, plainpages]{hyperref}
\usepackage{graphicx}

\usepackage{graphicx}
\usepackage{epsfig}

 \newtheorem{thm}{Theorem}[section]
 \newtheorem{cor}[thm]{Corollary}

 \newtheorem{prop}[thm]{Proposition}
 \theoremstyle{definition}
 
 \theoremstyle{remark}
 
 \numberwithin{equation}{section}
 
 \newcommand{\A}{\mathcal{A}}

\begin{document}

\title[Generalization  of Weak Amenability of group Algebras]
 {Generalization  of Weak Amenability of group Algebras}
\author{M. Eshaghi Gordji}
\address{
$^1$Department of Mathematics, Semnan University,\newline
\indent P. O. Box 35195-363, Semnan, Iran\newline
\indent Center of Excellence in Nonlinear Analysis and Applications (CENAA),\newline
\indent Semnan University, Iran}
\email{madjid.eshaghi@gmail.com}

\author{A. Jabbari}
\address{
$^2$ Young Researchers Club, Ardabil Branch\newline
\indent Islamic Azad University, Ardabil, Iran}
\email{jabbari\underline{ }al@yahoo.com}

\thanks{}

\thanks{}



\date{}

\dedicatory{}

\commby{}

\begin{abstract}
New notion of weak amenability, $(\varphi,\psi)$-weak amenability recently introduced. In this paper we consider this new notion and study $(\varphi,\psi)$-weak amenability of group algebras $M(G)$, $L^1(G)$ and $S^1(G)$.
\end{abstract}
\maketitle

\section{Introduction}

Suppose that $X$ is a Banach ${A}$-bimodule. A derivation
$D:{A}\rightarrow X$ is a linear map which satisfies
$D(ab)=a.D(b)+D(a).b$ for all $a,b\in {A}$.

 The derivation $\delta$
is said to be inner if there exists $x\in X$ such that
$\delta(a)=\delta_x(a)=a.x-x.a$ for all $a\in {A}$. The
linear space of bounded derivations from ${A}$ into $X$ denoted by
$Z^1({A},X)$ and the linear subspace of inner derivations denoted by
$N^1({A},X)$. We consider the quotient space
$H^1({A},X)=Z^1({A},X)/N^1({A},X)$, called the first Hochschild
cohomology group of ${A}$ with coefficients in $X$.

A Banach algebra $A$  is
amenable if every derivation from $A$ into every dual Banach
$A$-module is inner, equivalently if $H^1(A,X^*)=\{0\}$ for every
Banach $A$-module $X$,  this definition was introduced by Johnson
in \cite{jo1}. One of the important things that Johnson proved is that the group algebra $L^1(G)$ is amenable if and only if $G$ is amenable, where $G$ is locally compact group.

A Banach algebra ${A}$ is called weakly amenable if, $H^1({A},{A}^*)=\{0\}$ (for more details see \cite{Ba} and \cite{gr}). Let $A$ be  Banach algebra and let $\varphi$ and $\psi$ be
continuous homomorphisms on $A$. We consider the following module
actions on $A$,
\[a\cdot x: =\varphi(a)x , \quad x\cdot a: =x\psi(a)
\quad (a,x\in A).\]

 We denote  the above
$A$-module by $A_{(\varphi,\psi)}$.

Let $X$ be an $A$-module. A bounded linear mapping $d:A \longrightarrow X $ is
called a {\it$(\varphi,\psi)$-derivation} if
\[d(ab)=d(a)\cdot\varphi(b)+\psi(a)\cdot d(b)\quad(a,b\in
A).\]

A bounded linear mapping $d:A \to X $ is called a {\it
$(\varphi,\psi)$-inner derivation} if there exists $x\in X$ such
that
\[d(a)=x\cdot\varphi(a)-\psi(a)\cdot x \quad(a\in A).\]

Throughout, unless otherwise stated, by a $(\varphi,\psi)$-derivation we means that a continuous $(\varphi,\psi)$-derivation. New notion of weak amenability of Banach algebras, $(\varphi,\psi)$-weak amenability, introduced in \cite{go}, where the authors  determined the relations between weak amenability and $(\varphi,\psi)$-weak amenability of  Banach
algebras.

Let $A$ be a Banach algebra and
 let $\varphi$ and $\psi$ be continuous homomorphisms on $A$. Then $A$ is called $(\varphi,\psi)$-weakly amenable if $H^1(A,(A_{(\varphi,\psi)})^*)=\{0\}$.

Let $A$ and $B$ be Banach algebras. Similarly to \cite{go}, we
denote by $Hom(A,B)$ the metric space of all bounded homomorphisms
from $ A $ into $B $, with the metric derived from the usual
linear operator norm $\|\cdot\|$ on $ \mathcal B(A,B)$ and denote
$Hom(A,A)$ by $Hom(A)$, where $ \mathcal B(A,B)$ is the set of
all bounded linear operators from $A$ into $B$.

 Let $A$ be a Banach algebra, $X$ be a Banach $A$-module and
 let $\varphi,\psi\in Hom(A)$. A derivation $D:A\to X $ is called
 approximately
 $(\varphi,\psi)$-inner
if there exists a net $(x_{\alpha})$ in $X$ such that, for all
$a\in A$,
$D(a)=\lim_{\alpha}x_{\alpha}\cdot\varphi(a)-\psi(a)\cdot
x_{\alpha}$ in norm.  A  Banach algebra $A$ is
 approximately $(\varphi,\psi)$-weakly amenable if every derivation  $D:A\to (A_{(\varphi,\psi)})^*$ is
approximately $(\varphi,\psi)$-inner.

In this work we study $(\phi,\psi)$-weak amenability of Banach algebras. In section 2, we prove when the group algebra $L^1(G)$ is a two-sided $G$-module then it is $(\varphi,\psi)$-amenable. Also we consider the measure algebra $M(G)$, if $M(G)$ is $(\varphi,\psi)$-weakly amenable, then $G$ is discrete and amenable. In section 3, we investigate  $(\varphi,\psi)$-weak amenability of Segal algebras and abstract Segal algebras. At first, we consider approximate $(\varphi,\psi)$-derivation from Segal algebra $S^1(G)$ into its dual, and after that we study commutative abstract Segal algebras.

\section{$(\varphi,\psi)$-weak Amenability}
In this section, by providing some new results in $(\phi,\psi)$-weak amenability of Banach algebras, we study $(\varphi,\psi)$-weak amenability of group algebra $L^1(G)$ where $\varphi,\psi\in Hom (L^1(G))$ and $G$ is locally compact group.
\begin{prop}\label{2.1}
Let $A$ be a Banach algebra, $B$ be a Banach $A$-bimodule, $\Phi:A\longrightarrow B$ be a multiplier with dense range, and $d:A\longrightarrow (B_{(\varphi,\psi)})^*$ be a non-zero $(\varphi,\psi)$-derivation where $\varphi,\psi\in Hom(A)$. Then $D:=\Phi^*\circ d$ is a non-zero $(\varphi,\psi)$-derivation.
\end{prop}
\begin{proof}
For every $x,y,z\in A$ we have
\begin{eqnarray}
 \nonumber
   \langle D(xy),z\rangle&=&\langle \Phi^*\circ d(xy),z\rangle =\langle  \Phi^*(d(x).\varphi(y)+\psi(x).d(y)),z\rangle \\
 \nonumber
   &=&  \langle d(x).\varphi(y)+\psi(x).d(y),\Phi(z)\rangle= \langle d(x),\varphi(y).\Phi(z)\rangle+\langle d(y),\Phi(z).\psi(x)\rangle\\
 \nonumber
   &=& \langle d(x),\Phi(\varphi(y)z)\rangle+\langle d(y),\Phi(z\psi(x))\rangle=\langle \Phi^*(d(x)),\varphi(y)z\rangle+\langle \Phi^*(d(y)),z\psi(x)\rangle\\
    \nonumber
   &=& \langle D(x).\varphi(y)+\psi(x).D(y),z\rangle.
\end{eqnarray}

Therefore $D$ is a $(\varphi,\psi)$-derivation. If $D=0$, then for every $x,y\in A$ we have $\langle D(x),y\rangle=0$. Then $\langle  \Phi^*\circ d(x),y\rangle=\langle d(x),\Phi(y)\rangle=0$. This means that $d(A)=0$, and so $d=0$.
\end{proof}
\begin{thm}\label{2.2}
Let $A$ be a Banach algebra, $B$ be a closed subalgebra, and $I$ be an closed ideal of $A$ such that $A=B\oplus I$. If $A$ is $(\varphi,\psi)$-weakly amenable, where $\varphi,\psi\in Hom(A)$ and $\varphi(B),\psi(B)\subseteq B$, then $B$ is $(\varphi,\psi)$-weakly amenable.
\end{thm}
\begin{proof}
Let $\pi:A\longrightarrow B$ be the natural projection from $A$ onto $B$. Then $\pi$ is a multiplier on $B$. Also, for every $a,b\in A$, there are $x,y\in B$ such that $a=x+I$ and $b=y+I$, then $\pi(ab)=xy=x\pi(b)=\pi(a)y$.

Suppose that $d:B\longrightarrow (B_{(\varphi,\psi)})^*$ is an arbitrary $(\varphi,\psi)$-derivation. From Proposition \ref{2.1}, $D=\pi^*\circ d:A\longrightarrow (A_{(\varphi,\psi)})^*$ is a $(\varphi,\psi)$-derivation. Since $A$ is $(\varphi,\psi)$-weakly amenable, then there exists an element $\xi\in A^*$ such that
\begin{equation}\label{}
    \nonumber
    D(x)=\xi.\varphi(x)-\psi(x).\xi\hspace{1cm}(x\in A).
\end{equation}

Set $\eta=\xi|_B$. Then
\begin{eqnarray}
\nonumber
  \langle d(x),y\rangle &=& \langle d(x), \pi(y)\rangle =\langle \pi^*\circ d(x), y\rangle=\langle D(x), y\rangle\\
   \nonumber
   &=& \langle \xi.\varphi(x)-\psi(x).\xi,y\rangle =\langle \xi,\varphi(x)y\rangle-\langle \xi,y\psi(x)\rangle\\
   \nonumber
   &=& \langle \eta,\varphi(x)y\rangle- \langle \eta,y\psi(x)\rangle=\langle \eta.\varphi(x)-\psi(x).\eta,y\rangle,
\end{eqnarray}
for every $x,y\in B$ (note that $B$ is a closed subalgebra of $A$). Then $d$ is a $(\varphi,\psi)$-inner derivation.
\end{proof}

Let $G$ be a nondiscrete locally compact group. Then $M(G)=M_d(G)\oplus M_c(G)$ is the direct sum decomposition of the measure algebra $M(G)$ into its closed subalgebra of discrete measure $M_d(G)$ and its closed ideal of continuous measures $M_c(G)$ (Theorem 19.20 of \cite{he}).
\begin{cor}
Let $M(G)$ be $(\varphi,\psi)$-weakly amenable, where $\varphi,\psi\in Hom(M(G))$, such that $\varphi(M_d(G)),\psi(M_d(G))\subseteq M_d(G)$. Then $M_d(G)$ is $(\varphi,\psi)$-weakly amenable.
\end{cor}

In \cite[Theorem 2.1]{de}, Dales, Ghahramani and Helemskii showed that $M(G)$ is weakly amenable if and only if $G$ is discrete if and only if there is no non-zero, continuous point derivation at a character of $M(G)$.

 Similarly, let $\varphi\in Hom(M(G))$, and is onto. If $M(G)$ is $(\varphi,\varphi)$-weakly amenable and   $0\neq\psi\in\Delta(M(G))$, then by Theorem 2.10 of \cite{go}, there are no non-zero point derivation at $\psi\circ\varphi$. This means that $G$ is discrete and $M(G)$ is weakly amenable. Converse of this statement is true, when $\varphi$ and $\psi$ are identity homomorphisms.  This shows that the space of this new notion is wider then weak amenable Banach algebras.

Let $X$ be a Banach space, a net $(m_\alpha)\subset X^*$ called converges \emph{weak}$^\sim$ to $m\in X^*$ if $m_\alpha\stackrel{w^*}{\longrightarrow} m$ and $\|m_\alpha\|\to\|m\|$, this notion introduced in \cite{l1}. In particular case, if $\mu\in M(G)$, suppose that $\nu\in L^\infty(G)^*$ is a norm preserving extension of $\mu$. Then there exists a net $(f_\gamma)\subset L^1(G)$ with $\|f_\gamma\|\leq\|\mu\|$, and $f_\gamma\stackrel{w^*}{\longrightarrow} \nu$. By passing to a suitable subnet we can write $\|f_\gamma\|\to\|\mu\|$, so we have $f_\gamma{\longrightarrow} \mu$ weak$^\sim$. Let $A$ be a Banach algebra, for operator $\varphi$, we say that $\varphi$ satisfies in  weak$^\sim$ condition if for every net $(m_\alpha)\subset A^*$ which converges weak$^\sim$ to $m\in A^*$, then $m_\alpha\stackrel{w^*}{\longrightarrow}\varphi(m)$.
\begin{thm}\label{de}
Let $G$ be a locally compact group, and $X$ be a $M(G)$-bimodule by module actions $\mu.x=\widetilde{\varphi}(\mu).x$ and $x.\mu=x.\widetilde{\psi}(\mu)$ where $\varphi, \psi\in Hom(L^1(G))$, and $\widetilde{\varphi}, \widetilde{\psi}$ are extensions of $\varphi$ and $\psi$. Then every $(\varphi,\psi)$-derivation $D:L^1(G)\longrightarrow X^*$ extends to a unique $(\widetilde{\varphi},\widetilde{\psi})$-derivation form $M(G)$ into $X^*$.
\end{thm}
\begin{proof}
Let $\mu\in M(G)$, take a net $(f_\gamma)\subset L^1(G)$ such that $f_\gamma\to\mu$ weak$^\sim$. Then the net $(Df_\gamma)$ converges in w$^*$-topology. Define $\overline{D}\mu=w^*-\lim_\gamma Df_\gamma$. Thus $\overline{D}$ is a bounded linear operator which extend $D$. For every $x\in X$, $g_1,g_2\in L^1(G)$ and $\mu\in M(G)$ we have
\begin{eqnarray}\label{11}
 \nonumber
  \langle \overline{D}(\mu), g_1.x.g_2\rangle &=& w^*-\lim_\gamma\langle D(f_\gamma),g_1.x.g_2\rangle=w^*-\lim_\gamma\langle D(f_\gamma).g_1,x.g_2\rangle \\
    \nonumber
   &=& w^*-\lim_\gamma\langle D(f_\gamma).\varphi(g_{1}),x.g_2\rangle \\
 \nonumber
   &=& w^*-\lim_\gamma\langle D(f_\gamma.g_{1})-\psi(f_\gamma).D(g_{1}),x.g_2\rangle \\
   \nonumber
   &=&  w^*-\lim_\gamma\langle D(f_\gamma.g_{1}),x.g_2\rangle- w^*-\lim_\gamma\langle\psi(f_\gamma).D(g_{1}),x.g_2\rangle\\
   \nonumber
   &=&\langle \overline{D}(\mu.g_1),x.g_2\rangle-\langle \widetilde{\psi}(\mu).D(g_1),x.g_2\rangle\\
   &=& \langle \overline{D}(\mu.g_1),x.g_2\rangle-\langle D(g_1),x.g_2.\widetilde{\psi}(\mu)\rangle.
\end{eqnarray}

Now, we show that $\overline{D}$ is a $(\widetilde{\varphi},\widetilde{\psi})$-derivation. As above suppose that $f_\gamma\to\mu$ weak$^\sim$, then by \ref{11}, we have
\begin{eqnarray}
\nonumber
  \langle \overline{D}(\eta\mu),g.x\rangle &=& \lim_\gamma\langle {\overline{D}}(\eta.f_\gamma),g.x\rangle=\lim_\gamma\langle {\overline{D}}(\eta),f_\gamma g.x\rangle +\lim_\gamma\langle D(f_\gamma), g.x.\widetilde{\psi}(\eta)\rangle\\
   \nonumber
   &=& \lim_\gamma\langle {\overline{D}}(\eta).\varphi(f_\gamma), g.x\rangle +\lim_\gamma\langle D(f_\gamma), g.x.\widetilde{\psi}(\eta)\rangle\\
   \nonumber
   &=& \lim_\gamma \langle {\overline{D}}(\eta).\varphi(f_\gamma), g.x\rangle +\lim_\gamma\langle \widetilde{\psi}(\eta).D(f_\gamma), g.x\rangle\\
   \nonumber
   &=& \langle {\overline{D}}(\eta).\varphi(\mu)+ \widetilde{\psi}(\eta).\overline{D}(\mu), g.x\rangle,
\end{eqnarray}
for every $\eta\in M(G), g\in L^1(G)$, and $x\in X$ (note that $\|f_\gamma*g-\mu*g\|\to0$ and $L^1(G).X=X$). For uniqueness of $\overline{D}$, let $D'$ be other $(\widetilde{\varphi},\widetilde{\psi})$-derivation which extend $D$ to $M(G)$. Let $f_\gamma$ and $\mu$ be as above, then $D'(\mu)=\lim_\gamma D(f_\gamma)=\overline{D}(\mu)$.
\end{proof}

Let $G$ be a locally compact group. A Banach space $X$ called $G$-module if the following statements hold
\begin{enumerate}
  \item there is $k\geq0$ such that $\|g.x\|\leq k\|x\|$ for every $g\in G, x\in X$;
  \item for $x\in X$, the map $G\to X: g\mapsto g.x$ is continuous.
\end{enumerate}

Similarly for right Banach $G$-modules, and two-sided Banach $G$-modules, where in the latter case we require the map $G\times G\to X: (g_1,g_2)\mapsto g_1.x.g_2$ to be continuous. If $X^*$ is the dual of $X$ then by following actions, $X^*$ is two-sided $G$-module:
\begin{equation}\label{}
    \nonumber
    \langle f.\theta, x\rangle=\langle f,\theta.x\rangle\hspace{0.5cm}\emph{\emph{and}}\hspace{0.5cm}\langle \theta.f,x\rangle=\langle f,x.\theta\rangle,
\end{equation}
for every $\theta\in G, x\in X$ and $f\in X^*$.

Weak amenability of $L^1(G)$ studied by Johnson in \cite{jo2}, and after he, Despi\'{c} and Ghahramani in \cite{de} (Theorem 1) gave different proof for weak amenability of $L^1(G)$. We use technique of proof in \cite{de}, and we have the following Theorem.
\begin{thm}
Let $L^1(G)$ be two-sided $G$-module, and $\varphi,\psi\in Hom(L^1(G))$. Then $L^1(G)$ is $(\varphi,\psi)$-weak amenable.
\end{thm}
\begin{proof}
Let $D:L^1(G)\longrightarrow L^\infty(G)$ be a $(\varphi,\psi)$-derivation. By Theorem \ref{de}, it suffices to show extending of $D$ on $M(G)$ is inner. For $t\in G$, by $\delta_t$, we mean that the point mass at $t$. Then
\begin{eqnarray}\label{2.4}
\nonumber
   \psi(\delta_{t^{-1}}).D(\delta_t)&=& \psi(\delta_{t^{-1}}).D(\delta_{tx^{-1}}*\delta_x) \\
   \nonumber
   &=& \psi(\delta_{t^{-1}})\psi(\delta_{tx^{-1}}).D(\delta_x)+\psi(\delta_{t^{-1}}).D(\delta_{tx}).\varphi(\delta_x)\\
   &=& \psi(\delta_{x^{-1}}).D(\delta_x)+\psi(\delta_{x^{-1}}).{\big{(}}(\psi(\delta_{(tx^{-1})^{-1}}).D(\delta_{tx}){\big{)}}).\varphi(\delta_x)
\end{eqnarray}

For every $\lambda\in L^\infty(G)$, let $Re(\lambda)$ denote the real part of $\lambda$, and let $$S=\{Re{\big{(}}(\psi(\delta_{(tx^{-1})^{-1}}).D(\delta_{tx}){\big{)}}): t\in G\}.$$

Take $\xi=\sup(S)$. Since $L^1(G)$ is two-sided $G$-module, then we have
\begin{eqnarray}\label{2.5}
\nonumber
   \sup(\psi(\delta_{x^{-1}}).S.\varphi(\delta_x))&=&  \psi(\delta_{x^{-1}}).\sup(S).\varphi(\delta_x),~~~~~~~~~~~\emph{\emph{and}}\\
   \sup(\lambda+S) &=&\lambda+\sup(S)  \hspace{0.5cm}(x\in G, \lambda\in L^\infty(G)).
\end{eqnarray}

Now by \ref{2.4} and \ref{2.5}, we have
\begin{equation}\label{}
    \nonumber
    \xi= \psi(\delta_{x^{-1}}).Re(D(\delta_x))+\psi(\delta_{x^{-1}}).\xi.\varphi(\delta_x),
\end{equation}
then
\begin{equation}\label{}
\nonumber
   Re(D(\delta_x))= \psi(\delta_{x}).\xi-\xi.\varphi(\delta_x),
\end{equation}
for every $x\in G$. Similarly for imaginary part, there exists an element $\zeta\in L^\infty(G)$ such that
\begin{equation}\label{}
\nonumber
   Im(D(\delta_x))= \psi(\delta_{x}).\zeta-\zeta.\varphi(\delta_x),
\end{equation}
for every $x\in G$. Therefore by taking $\varsigma=\xi+i\zeta$ we have
\begin{equation}\label{}
\nonumber
   D(\delta_x)= \psi(\delta_{x}).\varsigma-\varsigma.\varphi(\delta_x),
\end{equation}
for every $x\in G$. Let $\mu\in M(G)$, then there exists a net $(\mu_\alpha)$ with each $\mu_\alpha$ a linear combination of point masses such that $\mu_\alpha\to\mu$ in strong topology (note that point masses are extreme point of $M(G)$). Also we can give this net such as net in proof of Theorem \ref{de}. Then
$$\overline{D}(\mu)=\widetilde{\psi}(\mu).\varsigma-\varsigma.\widetilde{\varphi}(\mu),$$
for every $\mu\in M(G)$. This means that $L^1(G)$ is $(\varphi,\psi)$-weak amenable.
\end{proof}
\section{Results For Segal Algebras}
Let $(\mathcal{A},\|.\|)$ be a Banach algebra. Then $(\mathfrak{B},\|.\|^{'})$ is an abstract Segal algebra
with respect to $(\mathcal{A},\|.\|)$ if:

$(1)$ $\mathfrak{B}$ is a dense left ideal in $\mathcal{A}$ and $\mathfrak{B}$ is a Banach algebra with respect to $\|.\|^{'}$;

$(2)$ There exists $M>0$ such that $\|b\|\leq M\|b\|^{'}$, for each $b\in\mathfrak{B}$;

$(3)$ There exists $C>0$ such that $\|ab\|^{'}\leq C\|a\|^{'}\|b\|^{'}$, for each $a, b\in\mathfrak{B}$.

 Let $G$ be a locally compact group. A linear subspace $S^1(G)$ of $L^1(G)$ is said to be
 a Segal algebra, if it satisfies the following conditions:

 $(i)$ $S^1(G)$ is a dense  in $L^1(G)$;

 $(ii)$ If $f\in S^1(G)$, then $L_xf\in S^1(G)$, i.e. $S^1(G)$ is left translation invariant;

 $(iii)$ $S^1(G)$ is a Banach space under some norm $\|.\|_S$ and $\|L_xf\|_s=\|f\|_s$, for all $f\in S^1(G)$ and $x\in G$;

 $(iv)$ Map $x\mapsto L_xf$ from $G$ into $S^1(G)$ is continuous.

For more details about Segal algebras see \cite{ri}. Dales and Pandey studied Weak amenability of special case of Segal algebras in \cite{da1}, and after they, Weak amenability of Segal algebras and Lebesgue-Fourier algebra of a locally compact group improved by Ghaharamani and Lau in \cite{gh1}.
\begin{thm}\label{1}
Let $G$ be an amenable group and let $S^1(G)$ be a symmetric Segal algebra with approximate identity $(e_\alpha)_\alpha$ contained in center of $L^1(G)$. Give $\varphi,\psi\in Hom(S^1(G))$ such that extensions of $\varphi, \psi$ on $L^1(G)$ are continuous, and $\psi$ is onto.  Then every  bounded $(\varphi,\psi)$-derivation from $S^1(G)$ into $S^1(G)^*$ is approximately $(\varphi,\psi)$-inner.
\end{thm}
\begin{proof}
Let $(e_\alpha)_\alpha$ be an approximate identity of $S^1(G)$ contained in the center of $L^1(G)$, and let $\widetilde{\varphi}$ and $\widetilde{\psi}$ be extensions of $\varphi$ and $\psi$ on $L^1(G)$.

 For each $\alpha$, define $D_\alpha:L^1(G)\longrightarrow S^1(G)^*$ by $D_\alpha(f)=D(e_\alpha *f)-D(e_\alpha)\cdot\widetilde{\varphi}(f)$, for every $f\in L^1(G)$. At first we show that $D_\alpha$ is a bounded $(\varphi,\psi)$-derivation. Boundedness of $D_\alpha$ comes from that $L^1(G)$ acts continuously on $S^1(G)$ on the right, and so $f\longmapsto D(e_\alpha*f)$ is continuous from $L^1(G)$ into $S^1(G)^*$. Similarly $L^1(G)$ acts continuously on $S^1(G)$ on the left, and so $f\longmapsto D(e_\alpha)\cdot\widetilde{\varphi}(f)$, therefore $D_\alpha$ is continuous.

  Now, we should show that $D_\alpha$ is a $(\varphi, \psi)$-derivation. Let $f_1,f_2\in L^1(G)$, then
\begin{eqnarray}
\nonumber
  D_\alpha(f_1*f_2) &=&  D(e_\alpha*f_1*f_2)-D(e_\alpha)\cdot\widetilde{\varphi}(f_1*f_2)\\
   \nonumber
   &=&  norm-\lim_\beta(D(e_\alpha*f_1*e_\beta*f_2)-D(e_\alpha)\cdot\widetilde{\varphi}(f_1*f_2))\\
   \nonumber
   &=& norm-\lim_\beta(D(e_\alpha*f_1)\cdot\varphi(e_\beta*f_2)+\psi(e_\alpha*f_1)\cdot D(e_\beta*f_2)\\
   \nonumber
   && \hspace{8cm}-D(e_\alpha)\cdot\widetilde{\varphi}(f_1*f_2)) \\
   \nonumber
   &=&  D(e_\alpha*f_1)\cdot\widetilde{\varphi}(f_2)-D(e_\alpha)\cdot\widetilde{\varphi}(f_1*f_2)+w^*-\lim_\beta \psi(f_1*e_\alpha)\cdot D(e_\beta*f_2)\\
    \nonumber
    &=& ( D(e_\alpha*f_1)-D(e_\alpha)\cdot\widetilde{\varphi}(f_1))\cdot\widetilde{\varphi}(f_2)+w^*-\lim_\beta (\widetilde{\psi}(f_1) \psi(e_\alpha)\cdot D(e_\beta*f_2))\\
     \nonumber
     &=& ( D(e_\alpha*f_1)-D(e_\alpha)\cdot\widetilde{\varphi}(f_1))\cdot\widetilde{\varphi}(f_2)+w^*-\lim_\beta (\widetilde{\psi}(f_1)\cdot D(e_\alpha*e_\beta*f_2)\\
      \nonumber
     &&\hspace{7cm}-\widetilde{\psi}(f_1)\cdot D(e_\alpha)\cdot\varphi(e_\beta*f_2))\\
      \nonumber
      &=& ( D(e_\alpha*f_1)-D(e_\alpha)\cdot\widetilde{\varphi}(f_1))\cdot\widetilde{\varphi}(f_2)+\widetilde{\psi}(f_1)\cdot D(e_\alpha*f_2)-\widetilde{\psi}(f_1)\cdot D(e_\alpha)\cdot\widetilde{\varphi}(f_2)\\
      \nonumber
      &=& D_\alpha(f_1)\cdot\widetilde{\varphi}(f_2)+\widetilde{\psi}(f_1)\cdot D(f_2).
\end{eqnarray}

Therefore $D_\alpha$ is a $(\varphi, \psi)$-derivation. By Johnson Theorem there exists  $\xi_\alpha$ in $S^1(G)^*$ such that
\begin{equation}\label{}
    \nonumber
    D_\alpha(f)=\xi_\alpha\cdot\widetilde{\varphi}(f)-\widetilde{\psi}(f)\cdot\xi_\alpha=D(e_\alpha*f)-D(e_\alpha)\cdot\widetilde{\varphi}(f),
\end{equation}
for every $f\in L^1(G)$, and we have
\begin{equation}\label{}
    \nonumber
    D(e_\alpha)\cdot \widetilde{\varphi}(f)=D(e_\alpha*f)-{\psi}(e_\alpha)\cdot D(f)\stackrel{w^*}{\longrightarrow}0,
\end{equation}
for every $f\in S^1(G)$. Then
\begin{equation}\label{}
    \nonumber
    D(f)=w^*-\lim_\alpha \xi_\alpha\cdot\widetilde{\varphi}(f)-\widetilde{\psi}(f)\cdot\xi_\alpha\hspace{2cm}(f\in S^1(G)).
\end{equation}

If we take $f=e_\alpha$, since $(e_\alpha)$ is in the center of $L^1(G)$, then $D(e_\alpha)=0$. There for every $f\in S^1(G)$, we have
\begin{equation}\label{}
    \nonumber
    D(f)=norm-\lim_\alpha \xi_\alpha\cdot\widetilde{\varphi}(f)-\widetilde{\psi}(f)\cdot\xi_\alpha,
\end{equation}
and so proof is complete.
\end{proof}
\begin{thm}
Let $\A$ be a commutative Banach algebra, and let $\mathfrak{B}$ be an abstract Segal algebra of $\A$ with approximate identity $(e_\alpha)$. Suppose that $\varphi\in Hom(\mathfrak{B})$, with continuous extension $\widetilde{\varphi}$ on $\A$. If $\A$ is $(\widetilde{\varphi},\widetilde{\varphi})$-weakly amenable, then $\mathfrak{B}$ is $(\varphi,\varphi)$-weakly amenable.
\end{thm}
\begin{proof}
Let $D:\mathfrak{B}\longrightarrow\mathfrak{B}^*$ be a continuous $(\varphi,\varphi)$-derivation. Similarly to proof of Theorem \ref{1}, define $D_\alpha:\A\longrightarrow\mathfrak{B}^*$ by $D_\alpha(f)=D(e_\alpha *f)-D(e_\alpha)\cdot\widetilde{\varphi}(f)$, for every $f\in A$. According to the proof of Theorem \ref{1}, $D_\alpha$ is a continuous $(\varphi,\varphi)$-derivation. Since $\A$ is $(\widetilde{\varphi},\widetilde{\varphi})$-weakly amenable and commutative, therefore $\mathfrak{B}^*$ is a symmetric $\A$-bimodule. These mean that $D_\alpha=0$, and this means that $D=0$.
\end{proof}


\end{document}